\documentclass[a4paper,11pt]{amsart}
\usepackage{amssymb,xypic,latexsym,graphicx}
\usepackage[matrix,arrow]{xy}
\newtheorem{theo}{Theorem}

\newtheorem{theocite}{Theorem} 
\newtheorem{question}{Question} 
\newtheorem{hope}{Hope} 
\newtheorem{cor}{Corollary} 
\newtheorem{lemme}[subsubsection]{Lemma}
\newtheorem{prop}[subsection]{Proposition}
 
\newtheorem{conj}{Conjecture}

\newcommand{\cucu}{\hfill $\Box$}

\newcommand{\HH}{{\rm H}}

\newcommand{\CH}{{\rm CH}}
\newcommand{\ZZ}{{\mathbb Z}}
\newcommand{\QQ}{{\mathbb Q}}
\newcommand{\NN}{{\mathbb N}}
\newcommand{\RR}{{\mathbb R}}
\newcommand{\C}{{\mathbb C}}
\newcommand{\PP}{{\mathbb P}}

\newcommand{\im}{{\rm Im}\,}

\newcommand{\codim}{{\rm codim}\,}

\newcommand{\coker}{{\rm coker}\,}
\newcommand{\anneau}{{\mathcal O}}

\newcommand{\X}{{\mathcal X}}
\newcommand{\Y}{{\mathcal Y}}

\begin{document}
\baselineskip 14.1pt
\title{Asymptotic bounds for Nori's connectivity theorem}
\author{Ania Otwinowska}
\date{}
\maketitle
\addtocounter{section}{-1}

\begin{abstract} 
Let $Y$ be a smooth complex projective variety. I study the 
cohomology of smooth families of hypersurfaces $\X\to B$ for 
$B\subset\PP\HH^0(Y,\anneau(d))$ a codimension $c$ subvariety. 
I give an asymptotically optimal bound on $c$ and $k$ when
$d\to\infty$ for the space $\HH^k(Y\times B,\X,\QQ)$ to vanish, 
thus extending the validity of Lefschetz Hyperplane section 
Theorem and Nori's Connectivity Theorem~\cite{nori}. 
Next, I construct in the limit case explicit families of higher 
Chow groups whose class does not vanish in $\HH^k(Y\times B,\X,\QQ)$.
Some of them are indecomposable. This suggests that in the limit 
case the space $\HH^k(Y\times B,\X,\QQ)$ should be spanned by higher 
Chow groups, containing \cite{nori} and \cite{O} as special cases
({\it cf}. conjecture~\ref{cok}).
\end{abstract}
 
\section{Introduction}

\subsection{Nori's connectivity theorem}
Let $Y$ be a smooth complex proper algebraic variety of 
dimension $N+1$ with an ample invertible sheaf ${\mathcal L}$ and
$\pi\colon {\mathcal X}\to\PP\HH^0(Y,{\mathcal L})$ 
the universal family of hypersurfaces of class $c_1({\mathcal L})$. 
For any $b\in\PP\HH^0(Y,{\mathcal L})$, let $X_b=\pi^{-1}(b)$ be the 
corresponding hypersurface of $Y$. 

The Lefschetz hyperplane section theorem asserts that for 
$k+1\in\{0,\dots,N\}$,
$$\HH^{k+1}(Y,X_b,\ZZ)=0.$$

Nori's connectivity theorem extends this result 
to generically defined cohomology classes of degree all the way up to
$2N=\dim _{\RR}X$. More precisely, let $T$ be a smooth variety and 
$\phi\colon T\to \PP\HH^0(Y,{\mathcal L})$ a differentiable map. 
Let $\X_T=\X\times_{\PP\HH^0(Y,{\mathcal L})}T$ be the universal 
hypersurface over $T$ and let $\Y_T=Y\times T$. 
By the Leray spectral sequence, the Lefschetz Hyperplane 
Section Theorem implies that for $k+1\in\{0,\dots,N\}$,
$$\HH^{k+1}(\Y_T,\X_T,\ZZ)=0.$$

With stronger assumptions one can say more. Assume that $T$ is a smooth 
algebraic variety and that $\phi\colon T\to \PP\HH^0(Y,{\mathcal L})$ is 
a smooth (hence dominant) morphism. Let $\anneau(1)$ be 
a fixed very ample invertible sheaf on $Y$ and assume
${\mathcal L}=\anneau(d)$ for some positive integer $d$. 

\begin{theocite}[Nori's connectivity theorem \cite{nori}]
For $d$ large enough and $k+1\in\{0,\dots,2N\}$,
\begin{eqnarray*}\label{nono}
\HH^{k+1}(\Y_T,\X_T,\QQ)=0.
\end{eqnarray*}
\end{theocite}

Note that Nori's connectivity theorem holds more generally for families 
of complete intersections in $Y$ (instead of only hypersurfaces). 
The more general setting is motivated by the construction of algebraic
cycles which are homologically but not algebraically equivalent 
to zero. Since in this paper I do not consider applications to 
algebraic equivalence, I restrict myself to the case of hypersurfaces, 
which involves all the ideas and simplifies notations.      

Despite its topological nature, the proof of Nori's connectivity 
theorem relies on Hodge theory ({\it cf}. \cite{nori}, section~1). 
Namely, let $\Omega^p_{\Y_T,\X_T}$ be the
coherent sheaf of $\anneau_{Y_T}$-modules defined by the 
short exact sequence
$$0\to\Omega^p_{\Y_T,\X_T}\to\Omega^p_{\Y_T}
   \to j_{T*}\Omega^p_{\X_T}\to 0,$$
where $j_T\colon \X_T\to\Y_T$ is the closed immersion. 
Then $\HH^{p+q+1}(\Y_T,\X_T,\C)=
{\mathbb H}^{p+q+1}(\Y_T,\Omega^{\bullet}_{\Y_T,\X_T})$. 
Following Nori, let 
$$G^{p}\HH^{p+q+1}(\Y_T,\X_T,\C)={\rm Im}\bigl(
{\mathbb H}^{p+q+1}(\Y_T,\Omega^{\bullet\geq p}_{\Y_T,\X_T})
\to\HH^{p+q+1}(\Y_T,\X_T,\C)\bigr)$$ 
denote the Hodge-Grothendick filtration. There is a natural inclusion
$$F^p\HH^{p+q+1}(\Y_T,\X_T,\C)\subset 
G^{p}\HH^{p+q+1}(\Y_T,\X_T,\C),$$
where $F^{p}\HH^{p+q+1}(\Y_T,\X_T,\C)$ denotes the Hodge filtration 
of the mixed Hodge structure on $\HH^{p+q+1}(\Y_T,\X_T)$. 
Hence if for some positive integer $k$ one has
$\HH^{q+1}(\Y_T,\Omega^p_{\Y_T,\X_T})=0$ for all
$p+q+1\leq k+1$ and $q\leq (k+1)/2$, then by
the Hodge--Fr\"ohlicher spectral sequence 
$G^{k+1-[(k+1)/2]}\HH^{k+1}(\Y_T,\X_T,\C)=0$, hence 
$F^{k+1-[(k+1)/2]}\HH^{k+1}(\Y_T,\X_T,\C)=0$ and a mixed Hodge structures 
argument implies $\HH^{k+1}(\Y_T,\X_T,\QQ)=0$.
Hence the connectivity theorem follows from the following more precise
Hogde-theoretical statement

\begin{theocite}[Nori] For $d$ large enough, $p+q+1\leq 2N$ and $q+1\leq N$,
\begin{eqnarray*}\label{nonobis}
\HH^{q+1}(\Y_T,\Omega^p_{\Y_T,\X_T})=0.
\end{eqnarray*}
\end{theocite}

\medskip

\subsection{Vanishing theorems}
The first goal of this paper is to investigate 
the minimal assumptions for Nori's connectivity theorem to extend. 
Natural questions are:
\begin{itemize}
\item how large does $d$ need to be?  
\item what happens if the map $\phi\colon T\to\PP\HH^0(Y,\anneau(d))$ 
fails to be dominant, i.e if one considers families of hypersurfaces 
which are not general, but close to being general?
\item can the above assertions be made more precise for a fixed integer 
$k$ (in Nori's connectivity theorem), and for fixed 
integers $p$ and $q$ (in the Hodge-theoretical version)? 
\end{itemize}

The first question was studied in Nagel~\cite{nagel} and 
Voisin~\cite{voisin}. 
In this paper I give the following asymptotic answer
to the above questions.

Let $T$ be a smooth algebraic variety with an algebraic morphism 
$\phi\colon T\to\PP\HH^0(Y,\anneau(d))$, such that $\X_T$ is smooth.
I assume that for all 
$t\in T$ the cokernel of the differential map
$d\phi(t)\colon T_{T,t}\to T_{\PP\HH^0(Y,\anneau(d)),\phi(t)}$ is
of dimension at most an integer $c$. 

\begin{theo}\label{vide}
For all $\varepsilon\in]0,1[$ and for all $p\in\NN$ and
for $q+1\in\{0,\dots ,N\}$, there is an integer $D$ depending only on 
$\varepsilon$, $p$, $Y$ and $\anneau(1)$, such that for all $d\geq D$,
$$\HH^{q+1}(\Y_T,\Omega^p_{\Y_T,\X_T})=0$$
provided that
$$c\leq (1-\varepsilon)\frac{d^{N-q}}{(N-q)!}.$$
\end{theo}

A straightforward adaptation of Nori's argument sketched above
shows that in the numerical situation of Theorem~\ref{vide} 
$G^{p}\HH^{p+q+1}(\Y_T,\X_T,\C)=0$, hence 
$F^{p}\HH^{p+q+1}(\Y_T,\X_T,\C)=0$, and the following topological 
statement holds.

\begin{cor}
For all $\varepsilon\in]0,1[$ and $k+1\in\{0,\dots ,2N\}$ 
there is an integer $D$ depending only on 
$\varepsilon$, $Y$ and $\anneau(1)$ such that for all $d\geq D$,
$$\HH^{k+1}(\Y_T,\X_T,\QQ)=0$$
provided that  
$$c\leq(1-\varepsilon)\frac{d^{N-[k/2]}}{(N-[k/2])!}.$$
\end{cor}

The following remarks show that the bounds given by Theorem~\ref{vide} 
and the Corollary are close to being asymptotically optimal:

\subsubsection{Case $k+1\in\{0,\dots ,N\}$}\label{vanN}
By the Lefschetz hyperplane section theorem and the 
Leray spectral sequence one has  
$\HH^{k+1}(\Y_T,\X_T,\QQ)=0$ with no assumption on $c$.
(similarily, by the holomorphic Leray spectral sequence one has 
$\HH^{q+1}(\Y_T,\Omega^p_{\Y_T,\X_T})=0$ for $p+q< N$ 
with no assumption on $c$).

\subsubsection{Case $k+1\in\{N+1,\dots ,2N\}$}\label{NN}
I show in the next section (see Theorem \ref{bornedec} below) that under 
a weak assumption on $Y$, the bound on $c$ is asymptotically optimal
for any couple of fixed integers $p$ and $q$ such that
$N/2\leq q\leq N-1$ and $p\geq q$. This implies that 
the bound given by the Corollary is asymptotically optimal. 

\subsubsection{Case $k+1\geq 2N+1$} One has
$\HH^{2N+1}(\Y_T,\X_T,\QQ)\neq 0$ for any positive integer $c$, 
$d>N$ and $T=\X_{B}$, where $B\subset\PP\HH^0(Y,\anneau(1))$ is any 
smooth codimension $c$ subvariety contained in the locus of 
smooth hypersurfaces. Indeed, the condition $d>N$ implies
$\HH^{0,N}(X_b,\QQ)\neq 0$ for all $b\in B$, hence
({\it cf}. \cite{nori}, section 0.2)  
the cohomology class of the diagonal embedding of $\X_{B}$ in 
$\X_T=\X_{B}\times_{B}\X_{B}$ is an element of $\HH^{2N}(\X_T,\QQ)$
which does not belong to $j_T^*\HH^{2N}(\Y_T,\QQ)$, hence defines 
a non-zero element of $\HH^{2N+1}(\Y_T,\X_T,\QQ)$.

\subsection{Some examples of algebraic classes}\label{examples}
The aim of this section is to give explicit constructions of
{\em algebraic} classes in $F^p\HH^{p+q+1}(\Y_T,\X_T,\C)$ for 
$q+1\in\{0,\dots N\}$, $p\geq q$ and
$c$ close to the bound given by theorem~\ref{vide}.

\subsubsection{The regulator map}
For any smooth algebraic variety $X$ and positive integers 
$p\geq q$ the higher Chow groups 
$\CH^p(X,p-q)$ are defined in \cite{bloch}. In this paper they
are always implicitely assumed to be  with $\QQ$ coefficients.
They come with the Deligne class map
$${\rm cl}_{p, p-q}({\mathcal D})\colon\CH^p(X,p-q)\to 
  \HH^{p+q}_{\mathcal D}(X,\QQ(p)),$$
where $\HH^{p+q}_{\mathcal D}(X,\QQ(p))$ denotes Deligne 
cohomology. There is a natural map
$$r_{p, p-q}\colon\HH^{p+q}_{\mathcal D}(X,\QQ(p))\to
  F^p\HH^{p+q}(X,\C),$$
where $F^{p}\HH^{p+q}(X,\C)$ denotes the Hodge filtration.  
Composing ${\rm cl}_{p, p-q}({\mathcal D})$ with $r_{p, p-q}$ 
one gets the regulator map
$${\rm reg}_{p,p-q}\colon\CH^p(X,p-q)\to F^p\HH^{p+q}(X,\C).$$
Note that this map is identically zero for $p>q$ whenever $X$ is proper.

Consider the case $X=\X_T$ (note that $\X_T$ is {\em not} proper). Let 
$$\HH^{p+q}(\X_T,\C)_v=\coker\bigl(j_T^*\colon
\HH^{p+q}(\Y_T,\C)\to\HH^{p+q}(\X_T,\C)\bigr)$$
denote the {\em vanishing} cohomology, let
$F^p\HH^{p+q}(\X_T,\C)_v=F^p\HH^{p+q}(\X_T,\C)\cap\HH^{p+q}(\X_T,\C)_v$
 and let
$$\overline{{\rm reg}}_{p,p-q}\colon\CH^p(\X_T,p-q)\to F^p\HH^{p+q}(\X_T,\C)_v
\subset F^p\HH^{p+q+1}(\Y_{T},\X_T,\C)$$
be the map deduced from ${\rm reg}_{p,p-q}$ by projection 
on $\HH^{p+q}(\X_T,\C)_v$. 

An element of
$F^p\HH^{p+q}(\X_T,\C)_v$ (or $F^p\HH^{p+q+1}(\Y_{T},\X_T,\C)$)
will be called {\em algebraic} if it lies in the image of 
$\overline{{\rm reg}}_{p,p-q}$.

\subsubsection{Chow groups} 
In the case $p=q$ higher Chow groups are Chow groups and the 
regulator is the class map.
I then write $\CH^p(\X_T)$ for $\CH^p(\X_T,0)$,
${\rm cl}_p$ for ${\rm reg}_{p,0}$ and  
$\overline{{\rm cl}}_p$ for $\overline{{\rm reg}}_{p,0}$

The following is a variation of my result in \cite{O}.

\begin{theo}\label{nl} 
Let $p$, $q$, $b$ and $c$ be positive integers such that
$p+q=N$, $q\in\{[(N+1)/2],\dots,N-1\}$, $b\geq 1$ and
$c\leq b\frac{d^{N-q}}{(N-q)!}$. Assume
$\phi(T)\subset\PP\HH^0(Y,\anneau(d))$ lies in the locus 
parametrizing smooth hypersurfaces. Then  
then there is an integer $D$
depending only on $b$, $Y$ and $\anneau(1)$ such that 
for all $d\geq D$ the following holds.

If $N$ is even and $p=q=N/2$ then the space  
$F^p\HH^{p+q}(\X_T,\C)_v$ is spanned via $\overline{{\rm cl}}_p$
by classes of flat families of algebraic cycles in $\X_T$ 
of codimension $p$ and degree at most $b$.

If $q>N/2$ then $F^p\HH^{p+q}(\X_T,\C)_v=0$.
\end{theo}

\subsubsection{Higher Chow groups} 
I now fulfill the promise made in section~\ref{NN}.

\begin{theo}\label{bornedec}
Let $p$ and $q$ be positive integers such that $p\geq q$
and $q\in\{[(N+1)/2],\dots,N-1\}$. 
Assume that $Y\subset\PP\HH^0(Y,\anneau(1))^{\vee}$ contains a 
codimension $q+1$ linear subvariety $V$ and let 
$T\subset\PP\HH^0(Y,\anneau(d))$
be the subvariety parametrising smooth hypersurfaces containing $V$. 
Then for $d\gg 0$ the variety $T$ is of codimension 
$c\sim_{d\to\infty}\frac{d^{N-q}}{(N-q)!}$, and there is a 
higher cycle $Z\in\CH^p(\X_{T},p-q)$ such that
$\overline{{\rm reg}}_{p,p-q}(Z)\neq 0$.
\end{theo}

\noindent
Hence in this situation the bound on $c$ in Theorem~\ref{vide} 
and its corollary are asymptotically optimal. 

Note that the assumption on $Y$ holds for $Y=\PP^{N+1}_{\C}$, which
should be regarded as the main example.

\subsubsection{Indecomposable Higher Chow groups}
For all $i\in\{0,\dots ,p\}$ and $j\in\{0,\dots ,p-q\}$
there are natural product maps
$$\psi_{i,j}\colon\CH^{p-i}(\X_T,p-q-j)
\otimes\CH^i(\X_T,j)\to\CH^p(\X_T,p-q).$$
I say that a higher cycle in $\CH^p(\X_T,p-q)$ is {\it decomposable}
if it belongs to the subgroup spanned by the images of the maps 
$\psi_{i,j}$ for $(i,j)\not\in\{(0,0),(p,p-q)\}$, and that it is
{\it indecomposable} otherwise. 

The proof of Theorem~\ref{bornedec} relies on an explicit construction
of higher cycles in $\CH^p(\X_{T},p-q)$. More precisely, since 
$\CH^1(\X_T,1)\simeq\HH^0(\X_T,\anneau^*_{\X_T})\simeq
\HH^0(T,\anneau^*_{T})$, iterating the above map $p-q$ times I 
obtain the map
$$\psi^{p-q}_{1,1}\colon\CH^{q}(\X_T)\otimes\HH^0(T,\anneau^*_{T})^{p-q}\to
  \CH^p(\X_T,p-q).$$
For $p>q$ the higher cycles in $\CH^p(\X_{T},p-q)$ I construct 
are in the image of $\psi^{p-q}_{1,1}$. Since $\psi^{p-q}_{1,1}$ 
commutes to the restriction
to a fiber, it follows that the restriction of these higher cycles 
to any fiber $X_t$, $t\in T$, is decomposable.

My result should be compared to Voisin's~\cite{voisin}, who studies 
the case $Y=\PP^{N+1}_{\C}$, $p=N$, $q=N-1$, $d=2N$ and 
$T\subset\PP\HH^0(Y,\anneau(d))$ an open subset (this implies $c=0$).
She constructs a family of higher cycles in $\CH^N(\X_T,1)$ that
limits the validity of Nori's connectivity theorem and whose restriction 
to a very general fiber $X_t$ for $t\in T$ is indecomposable. This suggests 
that for all $p,\ q$ such that $p+q\geq N$ and $p\geq q$ there might
be higher cycles in $\CH^p(\X_{T},p-q)$ for 
$c\sim_{d\to\infty}\frac{d^{N-q}}{(N-q)!}$ whose restriction to a very 
general fiber $X_t$, $t\in T$, is indecomposable. Supporting this 
intuition, I prove the following slightly weaker result for 
$Y=\PP^{N+1}$, $N$ even, $q=N/2$ and $p=q+1$.

\begin{theo}\label{borneimpaire}
Assume $N$ even, $q=N/2$ and that 
$Y\subset\PP\HH^0(Y,\anneau(1))^{\vee}$ contains a codimension 
$q$ linear subvariety. 
Then for $d\gg 0$ there is a subvariety 
$U\subset\PP\HH^0(Y,\anneau(d))$ of 
codimension $c\sim_{d\to\infty}2\frac{d^{N-q}}{(N-q)!}$, and a higher 
cycle $Z_U\in\CH^{q+1}(\X_{U},1)$ such that
\begin{itemize}
\item $\overline{{\rm reg}}_{q+1,1}(Z_U)\neq 0$, and
\item the restriction of $Z_U$ to a very general fiber $X_u$ for $u\in U$ 
is indecomposable.
\end{itemize}
\end{theo} 

\noindent
The cycle $Z_U$ is constructed explicitely and
the proof of indecomposability relies on Theorem~\ref{nl}. 

I believe that Theorem~\ref{borneimpaire} should generalize as follows.

\begin{hope}  \label{high}
Let $p$ and $q$ be positive integers such that $p>q$ and  
$q\in\{[(N+1)/2],\dots,N-1\}$. Assume $Y=\PP^{N+1}$. 
Then for $d\gg 0$ there is a subvariety 
$U\subset\PP\HH^0(Y,\anneau(d))$ of 
codimension $c\sim_{d\to\infty}(p-q+1)\frac{d^{N-q}}{(N-q)!}$
and a higher cycle $Z_{U}\in\CH^{p}(\X_{U},p-q)$ such that
\begin{itemize}
\item $\overline{{\rm reg}}_{p,p-q}(Z_U)\neq 0$, and
\item the restriction of $Z_U$ to a very general fiber $X_u$, 
$u\in U$, is indecomposable.
\end{itemize}
\end{hope}

\noindent
Theorem~\ref{borneimpaire} is the case $p=N/2+1$, $q=N/2$.

Note that $c\sim_{d\to\infty}2\frac{d^{N-q}}{(N-q)!}$ in 
Theorem~\ref{borneimpaire} (and 
$c\sim_{d\to\infty}(p-q+1)\frac{d^{N-q}}{(N-q)!}$
in its conjectural generalisation) while
Theorem~\ref{bornedec} gives examples of decomposable
cycles in $\CH^p(\X_{T},p-q)$ with a non zero vanishing 
cohomology class for $c\sim_{d\to\infty}\frac{d^{N-q}}{(N-q)!}$. 
This suggests that the bound on $c$ in Theorem~\ref{borneimpaire} 
may not be optimal.

\begin{question} Is there an indecomposable higher cycle
$Z_U\in\CH^p(\X_{U},p-q)$ such that $\overline{{\rm reg}}_{p,p-q}(Z_U)\neq 0$ 
for $c\sim_{d\to\infty}\frac{d^{N-q}}{(N-q)!}$?
\end{question}

\subsection{Algebraicity conjecture}
The results of Section~\ref{examples} suggest that
for $c$ slightly bigger than the bound 
$\displaystyle{\frac{d^{N-q}}{(N-q)!}}$
given by Theorem~\ref{vide}
the space $F^p\HH^{p+q}(\X_T,\QQ)_v$ should be 
algebraic. 

More precisely, I believe the following should hold 

\begin{conj}\label{copq}
For all integers $b\in\NN^*$ and $p\in\NN$ there is an integer $D$ 
depending only on $b$, $p$, $Y$ and on $\anneau(1)$, such that for all 
$d\geq D$ and $0\leq q\leq N-1$, if
$$c\leq b\frac{d^{N-q}}{(N-q)!}$$
then the map 
$\overline{{\rm reg}}_{p,p-q}:\CH^p(\X_T,p-q)\otimes\C\to 
F^p\HH^{p+q}(\X_T,\C)_v$ is surjective. 

By convention, $\CH^p(\X_T,p-q)=0$ for $p<q$.
\end{conj}

By Theorem~\ref{nl}, Conjecture~\ref{copq} holds for $p+q=N$, $q\geq N/2$
provided that $\phi(T)\subset\PP\HH^0(Y,\anneau(d))$ lies 
in the locus parametrizing smooth hypersurfaces.

\medskip

Mimicking Nori's argument, it is easy to see that 
Conjecture~\ref{copq} implies the following purely topological 
statement.

\begin{conj}\label{cok}
For every positive integer $b$ there is an integer $D$
depending only on $b$, $Y$ and $\anneau(1)$ such that 
for all $d\geq D$ and $k\in\{N,\dots ,2N-1\}$ and for
$$c\leq b\frac{d^{N-[k/2]}}{(N-[k/2])!},$$
\begin{itemize}
\item if $k$ is even then the map
$\overline{{\rm cl}}_{k/2}\colon\CH^{k/2}(\X_T)\otimes\C\to 
F^{k/2}\HH^{k}(\X_T,\C)_v$ is surjective;
\item if $k$ is odd then the maps
\begin{eqnarray*}
\overline{{\rm reg}}_{(k+1)/2,1}\colon
\CH^{(k+1)/2}(\X_T,1)\otimes\C&\to&\HH^k(\X_T,\C)_v,\quad{\it and}\\
\psi\circ\overline{{\rm reg}}_{(k+1)/2,1}\colon
\CH^{(k+1)/2}(\overline{\X_T},1)\otimes\C&\to&
\HH^k(\overline{\X_T},\C)_v\to\HH^k(\X_T,\C)_v
\end{eqnarray*}
generate $\HH^{k+1}(\Y_T,\X_T,\C)_v$.

Here $\overline{\X_T}$ is the variety isomorphic to $\X_T$ as a 
real manifold and endowed with the opposite complex structure and 
$\psi$ is the geometric frobenus.
\end{itemize}
\end{conj}

\noindent
Note that for $p+q<N$ Conjecture~\ref{copq} and for
$k<N$ Conjecture~\ref{cok} are trivially true since the 
target cohomology spaces vanish ({\it cf}. section~\ref{vanN}).

For $k$ even Theorem~\ref{vide} implies that in the numerical
situation of Conjecture~\ref{cok}, one has
$F^{k/2+1}\HH^{k}(\Y_T,\X_T,\C)=0$, hence by a mixed Hodge 
structures argument the space $\HH^{k}(\Y_T,\X_T,\C)$ is of pure 
Hodge type $(k/2,k/2)$ (like in \cite{nori}, section 0.2), hence
the subspace $\HH^k(\X_T,\C)_v$ is a sub-Hodge structure of pure 
Hodge type $(k/2,k/2)$. Hence for $k$ even Conjecture~\ref{cok} is 
implied by Hodge Conjecture.

\begin{figure}
\centering{
\parbox[c]{14.5cm}{\centering{
\includegraphics[height=6cm]{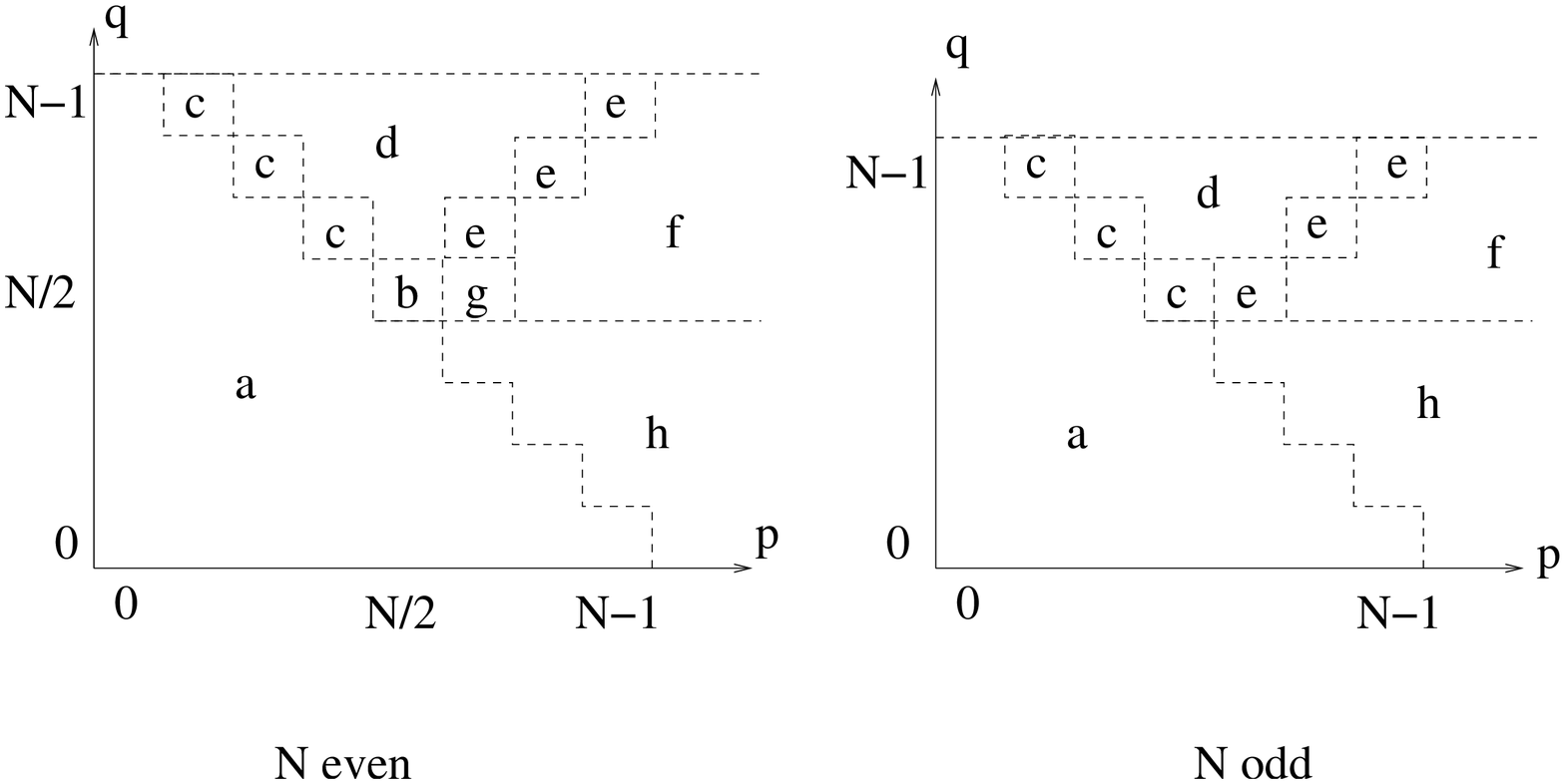} }} }
\end{figure}

\subsection{Summary} The figures summarise the results. 
Namely, for $c\leq b\frac{d^{N-q}}{(N-q)!}$

\begin{trivlist}
\item \ \  a~: $F^p\HH^{p+q+1}(\Y_T,\X_T,\C)=0$ 
  for any $c$ by section~\ref{vanN}. 
\item \ \  b~: $F^p\HH^{p+q+1}(\Y_T,\X_T,\C)=0$ 
  for $c\leq(1-\varepsilon){}{\frac{d^{N-q}}{(N-q)!}}$ by
  Theorem~\ref{vide} and $F^p\HH^{p+q}(\X_T,\C)_v$
  is spanned by classes of flat families of
  algebraic cycles of codimension $N/2$ and degree at most $b$ for
  $c\leq b{}{\frac{d^{N-q}}{(N-q)!}}$ by Theorem~\ref{nl},
  as predicted by Conjecture~\ref{copq}.
\item \ \  c~: $F^p\HH^{p+q+1}(\Y_T,\X_T,\C)=0$ 
  for $c\leq b{}{\frac{d^{N-q}}{(N-q)!}}$ by
  Theorem~\ref{nl}. 
\item \ \  d~: $F^p\HH^{p+q+1}(\Y_T,\X_T,\C)=0$ for 
  $c\leq(1-\varepsilon){}{\frac{d^{N-q}}{(N-q)!}}$ by
  Theorem~\ref{vide} and one should still have
  $F^p\HH^{p+q}(\X_T,\C)_v=0$ for 
  $c\leq b{}{\frac{d^{N-q}}{(N-q)!}}$ by
  Conjecture~\ref{copq}.
\item \ \  e~: $F^p\HH^{p+q+1}(\Y_T,\X_T,\C)=0$ for 
  $c\leq(1-\varepsilon){}{\frac{d^{N-q}}{(N-q)!}}$ by
  Theorem~\ref{vide} and $F^p\HH^{p+q}(\X_T,\C)_v$
  should be spanned by $\overline{\rm cl}_p(\CH^p(\X_T))$ for 
  $c\leq b{}{\frac{d^{N-q}}{(N-q)!}}$ by
  Conjecture~\ref{copq}; by Theorem~\ref{bornedec} there is a
  variety $Y$ such that $\overline{\rm cl}_p(\CH^p(\X_T))\neq 0$
  for $c\sim{}{\frac{d^{N-q}}{(N-q)!}}$.
\item \ \  f~: $F^p\HH^{p+q+1}(\Y_T,\X_T,\C)=0$ for $c\leq
  (1-\varepsilon){}{\frac{d^{N-q}}{(N-q)!}}$ by
  Theorem~\ref{vide} and $F^p\HH^{p+q}(\X_T,\C)_v$ should be spanned by 
  $\overline{\rm reg}_{p,p-q}\CH^p(\X_T,p-q)$ for 
  $c\leq b{}{\frac{d^{N-q}}{(N-q)!}}$ by Conjecture~\ref{copq}; 
  by Theorem~\ref{bornedec} there is a variety $Y$ such that 
  $\overline{\rm reg}_{p,p-q}\CH^p(\X_T,p-q)\neq 0$
  for $c\sim{}{\frac{d^{N-q}}{(N-q)!}}$ and by Conjecture~\ref{high}
  there should be an indecomposable higher cycle in 
  $\CH^p(\X_T,p-q)$ with a non zero image under 
  $\overline{\rm reg}_{p,p-q}$ for $c\sim(p-q+1){\frac{d^{N-q}}{(N-q)!}}$.
\item \ \  g~: as in f;
  by Theorem~\ref{borneimpaire} there is a variety $Y$ and 
  an indecomposable higher cycle in $\CH^p(\X_T,p-q)$ with 
  a non zero image under 
  $\overline{\rm reg}_{p,p-q}$ for $c\sim 2{}{\frac{d^{N-q}}{(N-q)!}}$. 
\item \ \  h~: $F^p\HH^{p+q+1}(\Y_T,\X_T,\C)=0$ 
  for $c\leq(1-\varepsilon){}{\frac{d^{N-q}}{(N-q)!}}$ by
  Theorem~\ref{vide} and $F^p\HH^{p+q}(\X_T,\C)_v$
  should be spanned by 
  $\overline{\rm reg}_{p,p-q}(\CH^p(\X_T,p-q))$ 
  for $c\leq b{}{\frac{d^{N-q}}{(N-q)!}}$ by Conjecture~\ref{copq}. 
  Here the situation is somewhat mysterious since I do not have 
  any example of a non-zero element in $\CH^p(\X_T,p-q)$ .
\end{trivlist}

\section{Proof of Theorem~\ref{vide} and of its Corollary} 

\subsection{Reduction to the local case}
A straightforward adaptation of Nori's argument shows that
if $\HH^{q+1}(Y_t,\Omega^p_{\Y_T,\X_T|Y_t})=0$ for all $t\in T$
then $\HH^{q+1}(\Y_T,\Omega^p_{\Y_T,\X_T})=0$ 
(\cite{nori}, Lemma 2.1). This 
shows that Theorem~\ref{vide} is implied by the following purely 
local-analytic statement.
 
\begin{prop}\label{vide'}
For all $\varepsilon\in\RR^*_+$ there is an integer $D$
depending only on $p$, $\varepsilon$, $Y$ and on $\anneau(1)$ such that 
\begin{eqnarray*}
  \HH^{q+1}(Y_t,\Omega^p_{\Y_T,\X_T|Y_t})=0
\end{eqnarray*}
for all integers $d\geq D$, $q\in\{-1,\dots ,N-1\}$, 
for all smooth analytic varieties $T$ endowed with a morphism 
$\phi\colon T\to\HH^0(Y,\anneau(d))$ such that 
\begin{itemize}
\item for all $t\in T$ the cokernel of $d\phi(t)$ is of dimension at 
most $(1-\varepsilon)\frac{d^{N-q}}{(N-q)!}$ and 
\item $\X_T$ is smooth. 
\end{itemize}
\end{prop}

\subsection{Proof of proposition~\ref{vide'}} 
Since the proof of proposition~\ref{vide'} 
proceeds by induction on several parameters, it is convenient to 
introduce the following notation. For a fixed positive integer $d$
and a projective space $\PP_{\C}^A$ I say that assertion 
${\mathcal H}(Y,p,q,c,a)$ holds for the smooth projective variety 
$Y\subset\PP_{\C}^A$ of dimension $N+1>0$ and for the integers 
$(p,q,c,a)\in\NN\times\{-1,\dots ,N-1\}\times\NN\times\ZZ$ if 
$$\HH^{q+1}(Y_t,\Omega^p_{\Y_T,\X_T|Y_t}\otimes\anneau(a'))=0$$
for all $a'\leq a$, for all smooth analytic varieties $T$ endowed 
with a morphism  $\phi\colon T\to\HH^0(Y,\anneau(d))$ such that $\X_T$ 
is smooth, and for all $t\in T$ such that the cokernel of $d\phi(t)$ 
is of dimension at most $c$.

The proof relies on the following four lemmas.

Recall ({\it cf}. \cite{greentheo}) that for any positive integer 
$c$, a $d$-decomposition of $c$ is the unique series 
$(c_d,\dots,c_r)$ such that $r\in \{ 1, \dots ,d\}$, 
$c_d>\cdots>c_r\geq r$ and
$$c= {c_d \choose d}+\cdots +{c_r \choose r}.$$
Let 
$$c_{<d>}={c_d-1 \choose d}+\cdots +{c_r-1 \choose r},$$ 
with the convention ${\alpha\choose\beta}=0$ for $\alpha<\beta$.

The asympotic behaviour of the function $c\mapsto c_{<d>}$ 
is described by the following lemma.
 
\begin{lemme}\label{greengeo}
For all integers $M\geq 2$ there exists a function 
$f_M\colon ]0,1[\times\NN\to ]0,1[$, $(\varepsilon,d)\mapsto\varepsilon'$ 
such that for all $\varepsilon\in]0,1[$ 
${\varepsilon'}\to{\varepsilon}^{\frac{M-1}{M}}$ for $d\to\infty$, 
and such that for all nonnegative integers 
$$c\leq(1-\varepsilon)\frac{d^{M}}{M!}$$
one has 
$$c_{<d>}\leq(1-\varepsilon')\frac{d^{M-1}}{(M-1)!}.$$
\end{lemme}

The following lemma is the key ingedient of my proof.
It relies on Green's hyperplane section theorem
\cite{greentheo}.

\begin{lemme}\label{recu}
Let $Y'$ be a generic hyperplane section of $Y$. 
Then for $q\leq N-2$
$$\left.\begin{array}{r}
     {\mathcal H}(Y',p,q,c_{<d>},a)\\
     {\mathcal H}(Y',p-1,q,c_{<d>},a-1)\end{array}\right\}
  \Rightarrow{\mathcal H}(Y,p,q,c,a).$$
By convention, assumption ${\mathcal H}(Y',-1,q,c_{<d>},a-1)$ holds.
\end{lemme}

The last two lemmas will only be used for $q=N-1$.
Both rely on an idea of Green~(\cite{koszul}, p. 193-194).

\begin{lemme}\label{c+1p-1} Assume $c>0$. Then
$$\left.\begin{array}{r}
     {\mathcal H}(Y,p+1,q,c-1,a)\\
     {\mathcal H}(Y,p+1,q-1,c,a)\end{array}\right\}
  \Rightarrow{\mathcal H}(Y,p,q,c,a).$$
By convention, assumption ${\mathcal H}(Y,p+1,-2,c,a)$ holds.
\end{lemme}

The last lemma has been proven by Nori for $a=0$ 
(\cite{nori}, remark 3.10); the general case involves no new ideas.
See also~\cite{nagel} for computations of an explicit bound for $s$.

\begin{lemme}\label{vide''} There is an integer $s$ depending only 
on $Y$ and $\anneau(1)$ such that for all $d\geq p+a+s$ the assertion 
${\mathcal H}(Y,p,q,0,a)$ holds for all $q\in\{-1,\dots,N-1\}$.
\end{lemme}

Lemma~\ref{vide''} will only be used for $q=N-1$.

\subsubsection{Proof of proposition~\ref{vide'}}
Postponing the proofs of the above Lemmas, 
I now show how they imply proposition~\ref{vide'}.

I show more generally that for all $\varepsilon\in]0,1[$
there is an integer $D$ depending only on $\varepsilon$, 
$p$, $a$ and $Y$ such that assertion ${\mathcal H}(Y,p,q,c,a)$ holds 
for all $d\geq D$, $q\in\{-1,\dots,N-1\}$
and $c\leq (1-\varepsilon)\frac{d^{N-q}}{(N-q)!}$.
The case $a=0$ is the proposition~\ref{vide'}.

I proceed by induction on $N-q$, $c$ and $N$.

Assume $N-q=1$. The condition on $c$ simplifies to 
$c\leq (1-\varepsilon)d$. Choose $D$ such that 
$\varepsilon D\geq p+a+s$. Thus for $d\geq D$ one has 
$d\geq (1-\varepsilon)d+\varepsilon D\geq c+p+a+s$.

The case $c=0$, $N-q=1$ follows from Lemma~\ref{vide''}. 

The case $c>0$, $N=0$  and $q=-1$ follows from Lemma~\ref{c+1p-1} by 
increasing induction on $c$. 
For the case $c>0$, $N-q=1$ let $Y'$ be a generic 
hyperplane section of $Y$. Lemmas~\ref{recu} and~\ref{c+1p-1}
provide successive implications
$$\left.\begin{array}{r}
     {\mathcal H}(Y,p+1,q,c-1,a)\\
\left.\begin{array}{r}
     {\mathcal H}(Y',p+1,q-1,c_{<d>},a)\\
     {\mathcal H}(Y',p,q-1,c_{<d>},a-1)\end{array}\right\}
\Rightarrow
     {\mathcal H}(Y,p+1,q-1,c,a)\end{array}\right\}
  \Rightarrow{\mathcal H}(Y,p,q,c,a).$$ 
Assertion ${\mathcal H}(Y,p,q,c,a)$ follows by increasing induction 
on both $c$ and $N$.  

Assume $N-q>1$ and let $Y'$ be a generic hyperplane section of $Y$. 
By induction on $N$ and by Lemma~\ref{greengeo} applied to 
$M=N-q$, assertions ${\mathcal H}(Y',p,q,c_{<d>},a)$ and
${\mathcal H}(Y',p-1,q,c_{<d>},a-1)$ hold.
Hence by Lemma~\ref{recu} assertion ${\mathcal H}(Y,p,q,c,a)$ holds.
\cucu

\subsection{Proof of the Lemmas}

\subsubsection{Proof of Lemma~\ref{greengeo}} 
Let $e$ be the smallest positive integer such that
$$(1-\varepsilon)\frac{d^{M}}{M!}\leq 
  {M+d \choose d}-{M+e \choose e}$$ 
for all $q\in\{0,\dots ,N-1\}$. One has 
$e\sim\varepsilon^{1/M}d$ for $d\to\infty$. 

The map $c\mapsto c_{<d>}$ is increasing. Thanks to the identities
\begin{eqnarray*}
{M+d \choose d}-{M+e \choose e}&=&
  \sum_{i=1}^{d-e} {M+d-i \choose d-i+1}\\
{M-1+d \choose d}-{M-1+e \choose d-e}&=&
  \sum_{i=1}^{d-e} {M-1+d-i \choose d-i+1}
\end{eqnarray*}
One gets
$$c_{<d>}\leq{M-1+d \choose d}-{M-1+e \choose e}.$$
The term on the right is asymptotically equivalent to 
$\varepsilon^{\frac{M-1}{M}}\frac{d^{M-1}}{(M-1)!}$ for $d\to\infty$. 
Define $\varepsilon'\in]0,1[$ by 
$\varepsilon'=\varepsilon^{\frac{M-1}{M}}$ if   
$\frac{d^{M-1}}{(M-1)!}\geq{M-1+d \choose d}-{M-1+e \choose e}$
and by the formula
$$(1-\varepsilon')\frac{d^{M-1}}{(M-1)!}=
  {M-1+d \choose d}-{M-1+e \choose e}$$
otherwise. Thus $\varepsilon'\in]0,1[$ and 
$\varepsilon'\sim\varepsilon^{\frac{M-1}{M}}$ for $d\to\infty$.
\cucu

\subsubsection{Proof of Lemma~\ref{recu}}
Let $i\colon Y'\to Y$ be the projective embedding, let $\Y'_T=Y'\times T$ and 
let $\X'_T=\X_T\times_{\Y_T}\Y'_T$. Fix $t\in T$. 

Since $Y'$ is generic, one may assume that $\X'_T$ is smooth and  
that $X_t$ meets $Y'_t$ transversally in $Y_t$. Restricting $T$ 
if necessary, one may further assume that for all $u\in T$, $X_u$ meets 
$Y'_u$ transversally in $Y_u$.
Hence there is a well-defined map $\phi'\colon T\to\PP\HH^0(Y',\anneau(d))$,
constructed as follows. 

Let $p_Y\colon\HH^0(Y,\anneau(d))\setminus\{0\}\to\PP\HH^0(Y,\anneau(d))$
be the natural projection and let
$$p_Y^*\phi\colon T\times_{\PP\HH^0(Y,\anneau(d))}\HH^0(Y,\anneau(d))
  \to\HH^0(Y,\anneau(d))$$ 
be the map deduced from $\phi$ by base change. Then $\phi'$ is 
defined by the map 
$$p_{Y'}^*\phi'\colon T\times_{\PP\HH^0(Y,\anneau(d))}\HH^0(Y,\anneau(d))
  \to\HH^0(Y',\anneau(d))$$ 
obtained by composing $p_Y^*\phi$ with the linear projection
$\HH^0(Y,\anneau(d))\to\HH^0(Y',\anneau(d))$
which maps to $0$ the polynomials vanishing on $Y'$. 

Let $E=p_Y^*\im(d\phi(t))$, {\it resp.} $E'=p_{Y'}^*\im(d\phi'(t))$. 
Let 
\begin{eqnarray*}
c&=\dim\coker(d\phi(t))&=\codim\left(E,\HH^0(Y,\anneau(d))\right),
  \quad {\it resp.} \\
c'&=\dim\coker(d\phi'(t))&=\codim\left(E',\HH^0(Y',\anneau(d))\right).
\end{eqnarray*}
Let $S=\bigoplus_{\delta\in\NN}S_{\delta}$, {\it resp.}  
$S'=\bigoplus_{\delta\in\NN}S'_{\delta}$ denote the homogenous 
polynomial ring of $\HH^0(Y,\anneau(1))$, {\it resp.} 
$\HH^0(Y',\anneau(1))$. Let $I\subset S$, {\it resp.}  
$I'\subset S'$ denote the homogenous ideal of $Y$, {\it resp.} $Y'$. 
One has $c=\codim(E+I_d,S_d)$, {\it resp.} 
$c'=\codim(E'+I'_d,S'_d)$. Since $Y'\subset Y$ is generic, 
Green's hyperplane theorem~\cite{greentheo} asserts that
$$c'\leq c_{<d>}.$$ 
Hence the cokernel of $d\phi'(t)$ is of dimension at most $c_{<d>}$.

The sheaf $\Omega^p_{\Y_T,\X_T|\Y'_T}$ fits into the short exact sequence
$$0\longrightarrow
  \Omega^{p-1}_{\Y'_T,\X'_T}\otimes\anneau_Y'(-1)
    \stackrel{dH\wedge}{\longrightarrow}
  \Omega^p_{\Y_T,\X_T|\Y'_T}\longrightarrow 
  \Omega^{p}_{\Y'_T,\X'_T}\longrightarrow 0$$ 
where $H\in\HH^0(Y,\anneau(1))$ is an equation of $Y'$.
Restricting this sequence to $Y'_t$ and tensorising it by 
$\anneau_Y'(a)$ one gets the associated long exact sequence 
\begin{eqnarray*}
\HH^{q+1}(Y'_t,\Omega^{p-1}_{\Y'_T,\X'_T|Y'_t}\otimes\anneau_Y'(a-1))
  \stackrel{dH\wedge}{\longrightarrow} 
\HH^{q+1}(Y'_t,\Omega^p_{\Y_T,\X_T|Y'_t}\otimes\anneau_Y'(a))
  \longrightarrow\\ \longrightarrow
\HH^{q+1}(Y'_t,\Omega^p_{\Y'_T,\X'_T|Y'_t}\otimes\anneau_Y'(a)).
\end{eqnarray*}
The term on the left vanishes by assumption 
${\mathcal H}(Y',p-1,q,c_{<d>},a-1)$ and the term on the right 
vanishes by assumption ${\mathcal H}(Y',p,q,c_{<d>},a)$.
Hence the middle term vanishes as well: one has
\begin{eqnarray}\label{recuaux}
\HH^{q+1}(Y'_t,\Omega^p_{\Y_T,\X_T|Y'_t}\otimes\anneau_Y'(a))=0.
\end{eqnarray}

For any coherent sheaf ${\mathcal F}$ on $Y_t$ there is a short 
exact sequence
$$0\longrightarrow{\mathcal F}\otimes\anneau(-1)
      \stackrel{H}{\longrightarrow}
  {\mathcal F}\longrightarrow i_*i^*{\mathcal F}\longrightarrow0,$$
Take ${\mathcal F}=\Omega^p_{\Y_T,\X_T|Y_t}\otimes\anneau_Y(a)$ and
consider the associated long exact sequence
\begin{multline*}
\HH^{q+1}(Y_t,\Omega^p_{\Y_T,\X_T|Y_t}\otimes\anneau_Y(a-1))
  \stackrel{H}{\longrightarrow} 
\HH^{q+1}(Y_t,\Omega^p_{\Y_T,\X_T|Y_t}\otimes\anneau_Y(a))\\
  \longrightarrow 
\HH^{q+1}(Y'_t,\Omega^p_{\Y_T,\X_T|Y'_t}\otimes\anneau_Y'(a)).
\end{multline*}
The term on the right vanishes by~(\ref{recuaux}). Since
$q\in\{0,\dots,N-1\}$, by increasing induction on $a$ the term on 
the right vanishes as well, the case $a\ll 0$ following from 
Serre's vanishing theorem. Hence the term in the middle vanishes.
\cucu

\subsubsection{Proof of Lemma~\ref{c+1p-1}}
Since the assertion depends only on the tangent space of $T$ at $t$,
one may replace $T$ by a small open neighbourhood of $t$ in the
tangent space of $T$ at $t$; hence one may assume $T=B\times S$,
where $B\subset\PP\HH^0(Y,\anneau(d))$ is an open 
neighbourghood of $\phi(t)$ of codimension at most $c$ and $S$
is an open ball.

Let $B'\subset\PP\HH^0(Y,\anneau(d))$ be a smooth analytic variety 
containing $B$ as a codimension~1 subvariety, and let 
$T'=B'\times S$. Then $\Y_{T}\subset \Y_{T'}$ is a smooth
codimension~1 subvariety and one has the following short exact 
sequence of sheaves on $\Y_T$:
$$0\longrightarrow
     \Omega^{p}_{\Y_{T},\X_{T}}\longrightarrow 
     \Omega^{p+1}_{\Y_{T'},\X_{T'}|\Y_{T}}\longrightarrow 
     \Omega^{p+1}_{\Y_{T},\X_{T}}\longrightarrow 0.$$ 
Restrict this sequence to $Y_t$, tensorise it with $\anneau(a)$ 
and consider the associated long exact sequence
$$\HH^{q}(Y_t,\Omega^{p+1}_{\Y_{T},\X_{T}|Y_t}\otimes\anneau(a))
  \longrightarrow 
  \HH^{q+1}(Y_t,\Omega^{p}_{\Y_{T},\X_{T}|Y_t}\otimes\anneau(a))
  \longrightarrow 
  \HH^{q+1}(Y_t,\Omega^{p+1}_{\Y_{T'},\X_{T'}|Y_t}\otimes\anneau(a)).$$ 
The terms on the right and on the left vanish by assumption, hence 
the term in the middle vanishes as well.
\cucu

\section{Proof of Theorem~\ref{nl}} 
Consider the exact sequence
$$R^k\pi_*\C\to R^k\pi_*(j_{T*}\C)\to R^k\pi_{*}(\C\to j_{T*}\C)\to 
  R^{k+1}\pi_*\C\to R^{k+1}\pi_*(j_{T*}\C)$$
and let ${\mathcal H}^{k}_{\X_T,v}=
\coker\left(R^k\pi_*\C\to R^k\pi_*(j_{T*}\C)\right)$
denote the locally trivial subsheaf of $R^k\pi_{Y*}(\C\to j_{T*}\C)$
whose fiber at any point $t\in T$ is 
canonically isomorphic to $\HH^{k}(X_t,\C)_v$.
Let $F^p\HH^0\left(T,{\mathcal H}^{k}_{\X_T,v}\right)$ denote the 
space of sections of ${\mathcal H}^{k}_{\X_T,v}$ 
whose fiber at any point $t\in T$ belongs to $F^p\HH^{k}(X_t,\C)_v$.

In \cite{O} I have shown

\begin{theocite}[\cite{O}]
For all integers $b\in\NN^*$ there is an integer $D$ depending only on 
$b$, $Y$ and on $\anneau(1)$, such that for all $d\geq D$, 
for all $q\in\{[(N+1)/2],\dots,N-1\}$, for all 
$$c\leq b\frac{d^{N-q}}{(N-q)!}$$
for all codimension $c$ algebraic varieties 
$T\subset\PP\HH^0(Y,\anneau(d))$
lying in the locus parametrizing smooth hypersurfaces
and for all $t\in T$, the following holds.

If $N$ is even and $q=N/2$ then the space 
$F^{N-q}\HH^0\left(T,{\mathcal H}^{N}_{\X_T,v}\right)$
is spanned by families of algebraic cycles of codimension $q$
and degree at most $b$.

If $q>N/2$ then
 $F^{N-q}\HH^0\left(T,{\mathcal H}^{N}_{\X_T,v}\right)=0$.
\end{theocite}

The proof of \cite{O} gives in fact a stronger result, namely
one can replace the subspace $T\subset\PP\HH^0(Y,\anneau(d))$ by 
any algebraic variety $T$ endowed with a morphism 
$\phi\colon T\to\PP\HH^0(Y,\anneau(d))$, such that 
\begin{itemize}
\item  $\phi(T)$ lies 
in the locus parametrizing smooth hypersurfaces and 
\item for all $t\in T$ the cokernel of the differential 
$d\phi(t)\colon T_{T,t}\to T_{\PP\HH^0(Y,\anneau(d)),\phi(t)}$ is
of dimension at most $c$;
\end{itemize}

This stronger version of the theorem of \cite{O} is
equivalent to Theorem~\ref{nl}. 
Indeed, one only needs to show
$$F^{N-q}\HH^{N}(\X_T,\C)_v\simeq
F^{N-q}\HH^0\left(T,{\mathcal H}^{N}_{\X_T,v}\right).$$
But the Lefschetz hyperplane Theorem and the Leray spectral 
sequence give a canonical isomorphism
$$\HH^{N}(\X_T,\C)_v\simeq 
\HH^0\left(T,{\mathcal H}^{N}_{\X_T,v}\right).$$
Since for any $t\in T$, the restriction map
$\HH^{N}(\X_T,\C)_v\to\HH^{N}(\X_t,\C)_v$
is a homomorphism of Hodge structures, hence is strict, the result 
follows. 
\cucu

\section{Proof of Theorem~\ref{bornedec}}\label{secbornedec}

Since $q\geq N/2$, the generic hypersurface of degree $d$ containing $V$
is smooth; let $T\subset\PP\HH^0(Y,\anneau(d))$ be a smooth open affine 
subvariety of the space of all smooth hypersurfaces containing $V$. 
One has $c\sim_{d\to\infty}\frac{d^{N-q}}{(N-q)!}$.

For all $i\in\{1,\dots,p-q\}$ choose hyperplane 
sections $H_i\subset T$ in general position and and let 
$f_i\in\HH^0(T,\anneau_{T})$ be the equations of $H_i$; Let 
$P_{i}=\bigcap_{j>i}H_j$, $Q_i=\bigcup_{j\leq i}H_j$ and 
$T_i=P_i\setminus\left(P_i\cap Q_i\right)$.  
For $j\leq i$ the functions $f_j$ do not vanish on $T_i$, hence belong
to $\HH^0(T_i,\anneau^*_{T_i})$ and there 
are well-defined higher cycles in $\CH^{q+i}(\X_{T_i},i)$
$$Z_i=\psi^{i}_{1,1}([V\times {T_i}],f_{1},\dots,f_{i}).$$
I show by increasing induction on $i$ that 
$\overline{{\rm reg}}_{q+i,i}(Z_i)\neq 0$. 
Case $i=p-q$ implies Theorem~\ref{bornedec}.

\medskip

Case $i=0$ a straightforward adaptation of \cite{voisin}, 
Proposition~4: for $d\gg 0$ and for any open subset 
$U_{0}\subset T_{0}$ the cohomology class of the cycle 
$[V\times U_{0}]\in\CH^q(\X_{U_{0}})$ 
in $\HH^{2q}(\X_{U_{0}},\C)$ does not belong to 
$j_{U_0}^*(\HH^{2q}(\Y_{U_{0}},\C))$. 

\medskip

Assume $i>0$. Then $T_{i-1}\cup T_i=
P _i\setminus\left(P_i\cap Q_{i-1}\right)$,
hence $T_{i-1}$ is a hyperplane section of $T_{i-1}\cup T_i$ and 
$T_i$ is the complementary open set. 
Hence there is a commutative diagram
$$\xymatrix{
\CH^{q+i}(\X_{T_i},i)\ar[r]^{l_{q+i,i}}\ar[d]^{\overline{{\rm reg}}_{q+i,i}}&
\CH^{q+i-1}(\X_{{T_{i-1}}},i-1)\ar[d]^{\overline{{\rm reg}}_{q+i-1,i-1}}\\
F^{q+i}\HH^{2q+i+1}(\Y_{T_{i}},\X_{T_{i}},\C)\ar[r]^{\rm res\quad}&
F^{q+i-1}\HH^{2q+i}(\Y_{T_{i}},\X_{T_{i-1}},\C)}$$
where res is the residue map and $l_{q+i,i}$ is the linking map for
the long  exact sequence of higher Chow groups. One has
$$l_{q+i,i}(Z_i)=
  (l_{q,0}\circ\psi^{i}_{1,1})([V\times {T_i}],f_{1},\dots,f_{i})=
  \psi^{i-1}_{1,1}([V\times {T_{i-1}}],f_{1},\dots,f_{i-1})=Z_{i-1}.$$
The induction assumption $\overline{{\rm reg}}_{q+i-1,i-1}(Z_{i-1})\neq 0$, and 
the commutativity of the diagram implies 
$\overline{{\rm reg}}_{q+i,i}(Z_i)\neq 0$. 
\cucu

\section{Proof of Theorem~\ref{borneimpaire}}

\subsection{Construction of the cycle $Z_U$}
Let $P\subset Y$ be a codimension~$q$ linear subspace, let 
$R\subset Y$ be a codimension~$q+2$ linear subspace such that
$R\subset P$ and let $S$ be the blow-up of $P$ along $R$. There 
are natural maps $\psi\colon S\to\PP^1_{\C}$ (where $\PP^1_{\C}$ parametrizes 
hyperplanes of $P$ containing $R$) and $s\colon S\to Y$. For any 
$x\in\PP^1_{\C}$ let $L_x=s(\psi^{-1}(x))$ denote the 
codimension~$q+1$ linear subspace of $Y$. For any subvariety 
$B\subset\PP\HH^0(Y,\anneau(d))$ let $C_B=S\times_Y\X_B$ and let 
$\psi_B\colon C_B\to\PP^1_{\C}$ and $s_B\colon C_B\to\X_B$ be the morphisms 
deduced from $\psi$ and $s$ by base change.

Let $U\subset\PP\HH^0(Y,\anneau(d))$ be the locus of hypersurfaces 
$X_t$ such that $[L_{0}\cap X_t]=dR$ and $[L_{\infty}\cap X_t]=dR$. 
One has $\codim(U,\PP\HH^0(Y,\anneau(d)))\sim_{d\to\infty}
2\frac{d^{N-q}}{(N-q)!}$. By construction, 
one has $s_U(\psi_U^{-1}(0))=R\times U=s_U(\psi_U^{-1}(\infty))$, 
hence the couple $(C_U,\psi_U)$ defines an element 
$Z_U\in\CH^{q+1}(U,1)$.

\subsection{Proof of $\overline{{\rm reg}}_{q+1,1}(Z_U)\neq 0$}
Let $T\subset\PP\HH^0(Y,\anneau(d))$ be the locus of hypersurfaces 
$X_t$ such that $L_{0}\subset X_t$ and $[L_{\infty}\cap X_t]=dR$.
Then $T\subset T\cup U$ is a hyperplane section and $U\subset T\cup U$ 
is the complementary open set. Hence there is a natural linking 
homomorphism $l_{q+1,1}\colon \CH^{q+1}(\X_{U},1)\to\CH^{q}(\X_T)$. 
Since $C_T$ is the union of a variety dominating $\PP^1_{\C}$ and of 
the variety $\psi_T^{-1}(0)\simeq L_0\times T$, one has
$l_{q+1,1}(Z_U)=d[s_T(\psi_T^{-1}(0))]=d[L_{0}\times T]$. 

I now proceed as in section~\ref{secbornedec}.
A straightforward adaptation of the argument of \cite{voisin}, 
Proposition 3 shows that for $d\gg 0$ the class 
${\rm cl}_q(d[L_{0}\times T])\in\HH^{2q}(\X_T,\C)$ does not belong to 
$j_T^*(\HH^{2q}(\Y_T,\C)$, hence $\overline{{\rm cl}}_q(d[L_{0}\times T])\in
\HH^{2q+1}(\Y_T,\X_T,\C)$ is non-zero.
Since the diagram
$$\xymatrix{\CH^{q+1}(\X_{U},1)\ar[r]^{l_{q+1,1}}\ar[d]^{\overline{{\rm reg}}_{q+1,1}}&
            \CH^{q}(\X_T)\ar[d]^{\overline{{\rm cl}}_q}\\
  \HH^{2q+2}(\Y_U,\X_U,\C)\ar[r]^{\rm res}&\HH^{2q+1}(\Y_T,\X_T,\C)}$$
is commutative, it follows that $\overline{{\rm reg}}_{q+1,1}(Z_U)\neq 0$.

\subsection{Indecomposability of the cycle $Z_U$}
Since $\overline{{\rm reg}}_{q+1,1}(Z_U)\neq 0$, it is enough to 
show that the space of decomposable cycles is contained in the kernel of
$\overline{{\rm reg}}_{q+1,1}$. 

Let $i\in\{1,\dots,q\}$, $Z'_U\in\CH^i(\X_U)$ and 
$Z''_U\in\CH^{q-i}(\X_U,1)$. One has to show 
${\rm reg}_{q+1,1}(\phi_{i,0}(Z'_U\otimes Z''_U))\in 
j_U^*F^{q+1}\HH^{2q+1}(\Y_U,\C)$.
Since the diagram 
$$\xymatrix{\CH^i(\X_U)\otimes\CH^{q+1-i}(\X_U,1)
          \ar[r]^{\qquad\phi_{i,0}}
          \ar[d]^{{\rm cl}_i\otimes{\rm reg}_{q+1-i,1}}&
        \CH^{q+1}(\X_U,1)\ar[d]^{{\rm reg}_{q+1,1}}\\
  F^i\HH^{2i}(\X_U,\C)\otimes F^{q-i+1}\HH^{2q-2i+1}(\X_U,\C)
          \ar[r]^{\qquad\qquad\quad\smile}&F^{q+1}\HH^{2q+1}(\X_U,\C),}$$
is commutative, it is enough to show ${\rm cl}_i(Z'_U)\in\HH^{2i}(\Y_U,\C)$ and 
${\rm reg}_{q+1-i,1}(Z''_U)\in\HH^{2q-2i+1}(\Y_U,\C)$.

If $i<q$ this follows from the isomorphisms given by Theorem~\ref{vide},
\begin{eqnarray*}
F^i\HH^{2i}(\X_U,\C)&\simeq &F^i\HH^{2i}(\Y_U,\C),\\
F^{q-i+1}\HH^{2q-2i+1}(\X_U,\C)&\simeq &F^{q-i+1}\HH^{2q-2i+1}(\Y_U,\C).
\end{eqnarray*}

If $i=q$, one has 
$$\CH^1(\X_U,1)\simeq\HH^0(U,\anneau^*(U))\simeq\CH^1(\Y_U,1),$$ 
hence ${\rm reg}_{1,1}(Z''_U)\in\HH^1(\Y_U,\C)$; on the other hand,
since $\X_U$ does not contain any flat family of cycles of 
codimension $q$ and degree less than two, by Theorem~\ref{nl} one has 
$F^q\HH^{2q}(\X_U,\C)\simeq F^q\HH^{2q}(\Y_U,\C)$, hence 
${\rm cl}_{q}(Z'_U)\in\HH^{2q}(\Y_U,\C)$.
\cucu

\subsection{Indecomposability of the restriction of $Z_U$ 
to a very general fiber $\X_u$, $u\in U$} This follows from the 
indecomposability of $Z_U$ by Proposition 5 of~\cite{voisin}. 

\bigskip

\noindent{\bf Aknowledgments}. I would like to thank Stefan M\"uller-Stach
for explaining to me what a higher Chow group is, and for a useful remark
on conjecture~\ref{copq}.

\end{document}